# A Proposed Artificial Neural Network Classifier to Identify Tumor Metastases

## Part I

M Khoshnevisan
Griffith University
Gold Coast, Queensland Australia

Sukanto Bhattacharya
Bond University
Gold Coast, Queensland Australia

Florentine Smarandache
University of New Mexico, USA

**Abstract:**

In this paper we propose a classification scheme to isolate truly benign tumors from those that initially start off as benign but subsequently show metastases. A non-parametric artificial neural network methodology has been chosen because of the analytical difficulties associated with extraction of closed-form stochastic-likelihood parameters given the extremely complicated and possibly non-linear behavior of the state variables. This is intended as the first of a three-part research output. In this paper, we have proposed and justified the computational schema. In the second part we shall set up a working model of our schema and pilot-test it with clinical data while in the concluding part we shall give an in-depth analysis of the numerical output and model findings and compare it to existing methods of tumor growth modeling and malignancy prediction.

**Introduction - mechanics of the mammalian cell cycle:**

The mammalian cell division cycle passes through four distinct phases with specific drivers, functions and critical checkpoints for each phase

| Phase | Main drivers | Functions | Checkpoints |
|---|---|---|---|
| $G_1$ (*gap* 1) | Cell size, protein content, nutrient level | Preparatory biochemical metamorphosis | Tumor-suppressor gene p53 |
| S (synthesization) | Replicator elements | New DNA synthesization | ATM gene (related to the MEC1 yeast gene) |
| $G_2$ (*gap* 2) | Cyclin B accumulation | Pre-mitosis preparatory changes | Levels of cyclin B/cdk1 – increased radiosensitivity |
| M (mitosis) | Mitosis Promoting Factor (MPF) – complex of cyclin B and cdk1 | Entry to mitosis; metaphase-anaphase transition; exit | Mitotic spindle – control of metaphase-anaphase transition |

The steady-state number of cells in a tissue is a function of the relative amount of cell proliferation and cell death. The principal determinant of cell proliferation is the residual effect of the interaction between oncogenes and tumor-suppressor genes. Cell death is determined by the residual effect of the interaction of proapoptotic and antiapoptotic genes. Therefore, the number of cells may increase due to either increased oncogenes activity or antiapoptotic genes activity or by decreased activity of the tumor-suppressor genes or the proapoptotic genes. This relationship may be shown as follows:

$$C_n = f(O, S, P, AP), \text{ such that } \{C_n'(O), C_n'(AP)\} > 0 \text{ and } \{C_n'(S), C_n'(P)\} < 0$$

Here $C_n$ is the steady-state number of cells, O is oncogenes activity, S is tumor-suppressor genes activity, P is proapoptotic genes activity and AP is antiapoptotic genes activity. The abnormal growth of tumor cells result from a combined effect of too few *cell-cycle decelerators* (tumor-suppressors) and too many *cell-cycle accelerators* (oncogenes). The most commonly mutated gene in human cancers is p53, which the cancerous tumors bring about either by overexpression of the p53 binding protein *mdm2* or through pathogens like the *human papilloma virus (HPV)*. Though not the objective of this paper, it could be an interesting and potentially rewarding epidemiological exercise to isolate the proportion of p53 mutation principally brought about by the overexpression of *mdm2* and the proportion of such mutation principally brought about by viral infection.

**Brief review of some existing mathematical models of cell population growth:**

Though the exact mechanism by which cancer kills a living body is not known till date, it nevertheless seems appropriate to link the severity of cancerous growth to the steady-state number of cells present, which again is a function of the number of oncogenes and tumor-suppressor genes. A number of mathematical models have been constructed studying tumor growth with respect to $C_n$, the simplest of which express $C_n$ as a function of time without any cell classification scheme based on histological differences. An **inherited cycle length model** was implemented by Lebowitz and Rubinow (1974) as an alternative to the simpler age-structured models in which variation in cell cycle times is attributed to occurrence of a chance event. In the LR model, variation in cell-cycle times is attributed to a distribution in inherited generation traits and the determination of the cell cycle length is therefore endogenous to the model. The population density function in the LR model is of the form **$C_n$ (a, t; τ)** where τ is the inherited cycle length. The boundary condition for the model is given as follows:

$$C_n (0, t; \tau) = 2 \int_0^\infty K(\tau, \tau') C_n (\tau', t; \tau') d\tau'$$

In the above equation, the kernel **K (τ, τ')** is referred to as the *transition probability function* and gives the probability that a parent cell of cycle length τ' produces a daughter cell of cycle length τ. It is the assumption that every dividing parent cell produces two daughters that yields the multiplier 2. The degree of correlation between the parent and daughter cells is ultimately decided by the choice of the kernel K. The LR model was further extended by Webb (1986) who chose to impose sufficiency conditions on the kernel K in order to ensure that the solutions asymptotically converge to a state of balanced exponential growth. He actually showed that the well-defined collection of mappings **{S (t): t ≥ 0}** from the Banach space B into itself forms a *strongly continuous semi-group of bounded linear operators*. Thus, for t ≥ 0, S (t) is the operator that transforms an initial distribution ϕ **(a, τ)** into the corresponding solution $C_n$ **(a, t; τ)** of the LR model at time t. Initially the model only allowed for a positive parent-daughter correlation in cycle times but keeping in tune with experimental evidence for such correlation possibly also being negative, a later; more general version of the Webb model has been developed which considers the sign of the correlation and allows for both cases.

There are also models that take $C_n$ as a function of both time as well as some physiological structure variables. Rubinow (1968) suggested one such scheme where the age variable "a" is replaced by a structure variable "μ" representing some physiological measure of cell maturity with a varying rate of change over time **v = dμ/dt**. If it is given that **$C_n$ (μ, t)** represents the cell population density at time t with respect to the structure variable μ, then the population balance model of Rubinow takes the following form:

$$\partial C_n/\partial t + \partial(vC_n)/\partial \mu = -\lambda C_n$$

Here λ **(μ)** is the maturity-dependent proportion of cells lost per unit of time due to non-mitotic causes. Either v depends on μ or on additional parameters like culture conditions.

**Purpose of the present paper:**
Growth in cell biology indicates changes in the size of a cell mass due to several interrelated causes the main ones among which are proliferation, differentiation and death. In a normal tissue, cell number remains constant because of a balance between proliferation, death and differentiation. In abnormal situations, increased steady-state cell number is attributable to either inhibited differentiation/death or increased proliferation with the other two properties remaining unchanged. Cancer can form along either route. Contrary to popular belief, cancer cells do not necessarily proliferate faster than the normal ones. Proliferation rates observed in well-differentiated tumors are not significantly higher from those seen in progenitor normal cells. Many normal cells hyperproliferate on occasions but otherwise retain their normal histological behavior. This is known as *hyperplasia*. In this paper, we propose a **non-parametric approach based on an artificial neural network classifier** to detect whether a hyperplasic cell proliferation could eventually become carcinogenic. That is, our model proposes to determine whether a tumor stays benign or subsequently undergoes metastases and becomes malignant as is rather prone to occur in certain forms of cancer.

**Benign versus malignant tumors:**

A benign tumor grows at a relatively slow rate, does not metastasize, bears histological resemblance to the cells of normal tissue, and tends to form a clearly defined mass. A malignant tumor consists of cancer cells that are highly irregular, grow at a much faster rate, and have a tendency to metastasize. Though benign tumors are usually not directly life threatening, some of the benign types do have the capability of becoming malignant. Therefore, viewed a stochastic process, a purely benign growth should approach some critical *steady-state mass* whereas any growth that subsequently becomes cancerous would fail to approach such a steady-state mass. One of the underlying premises of our model then is that cell population growth takes place according to the basic *Markov chain rule* such that the observed tumor mass in time $t_{j+1}$ is dependent on the mass in time $t_j$.

**Non-linear cellular biorhythms and chaos:**

A major drawback of using a parametric stochastic-likelihood modeling approach is that often closed-form solutions become analytically impossible to obtain. The axiomatic approach involves deriving analytical solutions of *stiff stochastic differential-difference* *equation systems*. But these are often hard to extract especially if the *governing system is decidedly non-linear* like Rubinow's suggested physiological structure model with velocity v depending on the population density $C_n$. The best course to take in such cases is one using a non-parametric approach like that of artificial neural networks.

The idea of chaos and non-linearity in biochemical processes is not new. Perhaps the most widely referred study in this respect is the *Belousov-Zhabotinsky (BZ) reaction*. This chemical reaction is named after B. P. Belousov who discovered it for the first time and A. M. Zhabotinsky who continued Belousov´s early work. R. J. Field, Endre Körös, and R. M. Noyes published the mechanism of this oscillating reaction in 1972. Their work opened an entire new research area of *nonlinear chemical dynamics*.

Classically the BZ reaction consist of a one-electron redox catalyst, an organic substrate that can be easily brominated and oxidized, and sodium or potassium bromate ion in form of $NaBrO_3$ or $KBrO_3$ all dissolved in sulfuric or nitric acid and mostly using $C_e$ (III)/$C_e$ (IV) salts and $M_n$ (II) salts as catalysts. Also Ruthenium complexes are now extensively studied, because of the reaction's extreme photosensitivity. There is no reason why the highly intricate intracellular biochemical processes, which are inherently of a much higher order of complexity in terms of molecular kinetics compared to the BZ reaction, could not be better viewed in this light. In fact, experimental studies investigating the physiological clock (of yeast) due to oscillating enzymatic breakdown of sugar, have revealed that the coupling to membrane transport could, under certain conditions, result in *chaotic biorhythms*. The yeast does provide a useful experimental model for histologists studying cancerous cell growth because the ATM gene, believed to be a critical checkpoint in the S stage of the cell cycle, is related to the MEC1 yeast gene. Zaguskin has further conjectured that all biorhythms have a *discrete fractal structure*.

The almost ubiquitous growth function used to model population dynamics has the following well-known *difference equation* form:

$$X_{t+1} = rX_t (1 - X_t/k)$$

Such models exhibit *period-doubling* and subsequently chaotic behavior for certain critical parameter values of r and k. The limit set becomes a fractal at the point where the model degenerates into pure chaos. We can easily deduce in a discrete form that the original Rubinow model is a linear one in the sense that $C_{nt+1}$ is *linearly dependent* on $C_{nt}$:

$$\Delta C_n/\Delta t + \Delta(vC_{nt})/\Delta\mu = -\lambda C_{nt}, \text{ that is}$$

$$(\Delta C_n/\Delta t) + (\Delta v/\Delta\mu) C_{nt} + (\Delta C_{nt}/\Delta\mu) v = -\lambda C_{nt}$$

$$\Delta C_n = -C_{nt}(\lambda + \Delta v/\Delta\mu) / (2/\Delta t) \ldots \text{ as } v = \Delta\mu/\Delta t, \text{ that is}$$

$$C_{nt+1} = rC_{nt}(1 - 1/k) \ldots \text{ putting } k = -[(2/\Delta t) - 1 - (\lambda + \Delta v/\Delta\mu)]^{-1} \text{ and } r = (2/\Delta t)^{-1}$$

Now this may be oversimplifying things and the true equation could indeed be analogous to the non-linear population growth model having a more recognizable form as follows:

$$C_{nt+1} = rC_{nt}(1 - C_{nt}/k)$$

Therefore, we take the conjectural position that very similar period-doubling limit cycles degenerating into chaos could explain some of the sudden "jumps" in cell population observed in malignancy when the standard linear models become drastically inadequate.

No linear classifier can identify a *chaotic attractor* if one is indeed operating as we surmise in the biochemical molecular dynamics of cell population growth. A non-linear and preferably non-parametric classifier is called for and for this very reason we have proposed artificial neural networks as a fundamental methodological building block here. Similar approach has paid off reasonably impressively in the case of complex systems modeling, especially with respect to weather forecasting and financial distress prediction.

**Artificial neural networks primer:**

Any artificial neural network is characterized by specifications on its *neurodynamics* and *architecture*. While neurodynamics refers to the input combinations, output generation, type of mapping function used and weighting schemes, architecture refers to the network configuration i.e. type and number of neuron interconnectivity and number of layers.

The input layer of an artificial neural network actually acts as a buffer for the inputs, as numeric data are transferred to the next layer. The output layer functions similarly except for the fact that the direction of dataflow is reversed. The transfer activation function is one that determines the output from the weighted inputs of a neuron by *mapping* the input data onto a suitable solution space. The output of neuron j after the summation of its weighted inputs from neuron 1 to i has been mapped by the transfer function f can be shown to be as follows:

$$O_j = f_j (\Sigma w_{ij} x_i)$$

A transfer function maps any real numbers into a domain normally bounded by 0 to 1 or −1 to 1. The most commonly used transfer functions are sigmoid, hypertan, and Gaussian.

A network is considered fully connected if the output from a neuron is connected to every other neuron in the next layer. A network may be *forward propagating* or *backward propagating* depending on whether outputs from one layer are passed unidirectionally to the succeeding or the preceding layer respectively. Networks working in closed loops are termed recurrent networks but the term is sometimes used interchangeably with backward propagating networks. Fully connected feed-forward networks are also called *multi-layer perceptrons* (*MLPs*) and as of now they are the most commonly used artificial neural network configuration. Our proposed artificial neural network classifier may also be conceptualized as a *recursive combination* of such *MLPs*.

Neural networks also come with something known as a hidden layer containing *hidden neurons* to deal with very complex, non-linear problems that cannot be resolved by merely the neurons in the input and output layers. There is no definite formula to determine the number of hidden layers required in a neural network set up. A useful heuristic approach would be to start with a small number of hidden layers with the numbers being allowed to increase gradually only if the learning is deemed inadequate. This should theoretically also address the regression problem of *over-fitting* i.e. the network performing very well with the training set data but poorly with the test set data. A neural network having no hidden layers at all basically becomes a *linear classifier* and is therefore statistically indistinguishable from the general linear regression model.

**Model premises:**

(1) The function governing the biochemical dynamics of cell population growth is inherently non-linear

(2) The sudden and rapid degeneration of a benign cell growth to a malignant one may be attributed to an underlying chaotic attractor

(3) Given adequate training data, a non-linear binary classification technique such as that of Artificial Neural Networks can *learn* to detect this underlying chaotic attractor and thereby prove useful in predicting whether a benign cell growth may subsequently turn cancerous

**Model structure:**

We propose a nested approach where we treat the output generated by an earlier phase as an input in a latter phase. This will ensure that the artificial neural network virtually acts as a *knowledge-based system* as it takes its own predictions in the preceding phases into consideration as input data and tries to generate further predictions in succeeding phases. This means that for a k-phase model, our set up will actually consist of k recursive networks having k phases such that the $j^{th}$ phase will have input function $\mathbf{I_j = f \{O (p'_{j-1}), I (p_{j-1}), p_j\}}$, where the terms $\mathbf{O (p'_{j-1})}$ and $\mathbf{I (p_{j-1})}$ are the output and input functions of the previous phase and $p_j$ is the vector of additional inputs for the $j^{th}$ stage. The said recursive approach will have the following schema for a nested artificial neural network model with k = 3:

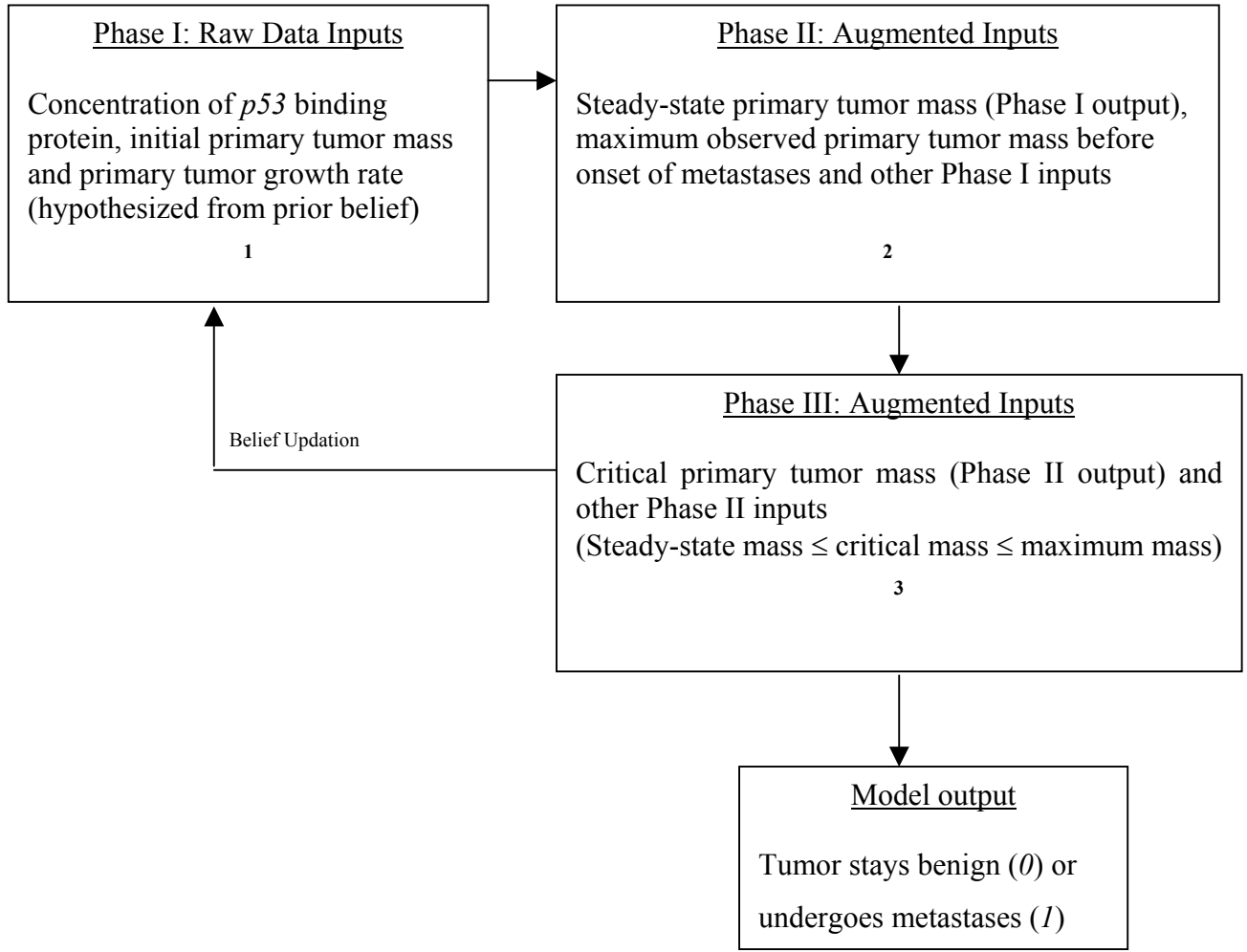

Phase I – target class variable: benign primary tumor mass

Phase II – target class variable: primary tumor mass at point of detection of malignancy

Phase III – target class variable: metastases ($M$) $\rightarrow$ *1*, no metastases ($B$) $\rightarrow$ *0*

As is apparent from the above schema, the model is intended to act as a sort of a *knowledge bank* that continuously keeps updating *prior* beliefs about tumor growth rate. The critical input variables are taken as concentration of p53 binding protein and observed tumor mass. The first one indicates the activity of the oncogenes vis-à-vis the tumor suppressors while the second one considers the extent of hyperplasia.

The model is proposed to be trained in phase I with histopathological data on concentration of p53 binding protein along with clinically observed data on tumor mass.

The inputs and output of Phase I is proposed to be fed as input to Phase II along with additional clinical data on maximum tumor mass. The output and inputs of Phase 2 is finally to be fed into Phase III to generate the model output – a binary variable *M|B* that takes value of *1* if the tumor is predicted to metastasize or *0* otherwise. The recursive structure of the model is intended to pick up any underlying chaotic attractor that might be at work at the point where benign hyperplasia starts to degenerate into cancer. Issues regarding network configuration, learning rate, weighting scheme and mapping function are left open to experimentation. It is logical to start with a small number of hidden neurons and subsequently increase the number if the system shows inadequate learning.

**Addressing the problem of training data unavailability:**

While training a neural network, if no target class data is available, the complimentary class must be inferred by default. Training a network only on one class of inputs, with no counter-examples, causes the network to classify everything as the only class it has been shown. However, by training the network on randomly selected counter-examples during training can make it behave as a *novelty detector* in the test set. It will then pick up any deviation from the norm as an *abnormality*. For example, in our proposed model, if the clinical data for initially benign tumors subsequently turning malignant is unavailable, the network can be trained with the benign cases with random inputs of the malignant type so that it automatically picks up any deviation from the norm as a possible malignant case.

A mathematical justification for synthesizing unavailable training data with random numbers can be derived from the fact that network training seeks to minimize the sum squared of errors over the training set. In a binary classification scheme like the one we are interested in, where a single input k produces an output f (k), the desired outputs are 0 if the input is a benign tumor that has stayed benign (*B*) and 1 if the input is a benign tumor that has subsequently turned malignant (*M*). If the prior probability of any piece of data being a member of class *B* is $P_B$ and that of class *M* is $P_M$; and if the probability distribution functions of the two classes expressed as functions of input k are $p_B(k)$ and $p_M(k)$, then the sum squared error, $\varepsilon$, over the entire training set will be given as follows:

$$\varepsilon = {}_{-\infty}\!\int^{\infty} P_B p_B (k)[f (k) - 0]^2 + P_M p_M (k)[f (k) -1]^2 \, dk$$

Differentiating this equation with respect to the function f and equating to zero we get:

$$\partial\varepsilon/\partial f = 2p_B (k) P_B f (k) + 2p_M (k) P_M f (k) - 2p_M (k) P_M = 0 \text{ i.e.}$$

$$f (k)^* = [p_M (k) P_M] / [p_B (k) P_B + p_M (k) P_M]$$

The above optimal value of f (k) is exactly the same as the probability of the correct classification being *M* given that the input was k. This shows that by training for minimization of sum squared error; and using as targets 0 for class *B* and 1 for class *M*, the output from the network converges to an identical value as the probability of class *M*.

**Gazing upon the road ahead:**

The main objective of our proposed model is to isolate truly benign tumors from those that initially start off as benign but subsequently show metastases. The non-parametric artificial neural network methodology has been chosen because of the analytical difficulties associated with extraction of closed-form stochastic likelihood parameters given the extremely complicated and possibly non-linear behavior of the state variables. This computational approach is proposed as a methodological alternative to the stochastic calculus techniques of tumor growth modeling commonly used in mathematical biology. Though how the approach actually performs with numerical data remains to be extensively tested, the proposed schema has been made as flexible as possible to suit most designed experiments to test its performance effectiveness and efficiency. In this paper we have just outlined a research approach – we shall test it out in a subsequent one.

*References:*

Atchara Sirimungkala, Horst-Dieter Försterling, and Richard J. Field, "Bromination Reactions Important in the Mechanism of the Belousov-Zhabotinsky Reaction", J. of Physical Chemistry, 1999.